\documentclass[a4paper,12pt]{amsart}

\usepackage{amssymb,amsbsy,amsmath,amsfonts,amssymb,amscd}
\usepackage{latexsym}
\usepackage{graphics}
\usepackage{color}
\input xy
\xyoption{all}

\newcommand\sC{{\mathcal C}}

\newcommand\sD{{\mathcal D}}

\newcommand\sI{{\mathcal I}}

\newcommand\s{\sigma}
\newcommand\Si{\Sigma}

\newcommand\Ga{\Gamma}

\newcommand{\CC}{\ensuremath{\mathbb{C}}}

\newcommand{\ZZ}{\ensuremath{\mathbb{Z}}}

\newcommand{\sS}{\ensuremath{\mathcal{S}}}

\newcommand{\hol}{\ensuremath{\mathcal{O}}}

\newcommand{\PP}{\ensuremath{\mathbb{P}}}

\newcommand{\ra}{\ensuremath{\rightarrow}}

\def\eea{\end{eqnarray*}}
\def\bea{\begin{eqnarray*}}

\newcommand\dual{\mathrel{\raise3pt\hbox{$\underline{\mathrm{\thinspace d
\thinspace}}$}}}
\newcommand\qe{\ifhmode\unskip\nobreak\fi\quad $\Box$}       

\def\BOX{\hfill\lower.5\baselineskip\hbox{$\Box$}}

\newtheorem{theo}{Theorem}[section]
\newtheorem{remarkk}[theo]{Remark}
\newtheorem{lemma}[theo]{Lemma}
\newtheorem{corollary}[theo]{Corollary}
\newenvironment{rem}{\begin{remarkk}\rm}{\end{remarkk}}

\newtheorem{defin}[theo]{Definition}

\newtheorem{example}[theo]{Example}

\newcommand{\Proof}{{\it Proof. }}

\begin{document}
\title[Uniformization by a symmetric domain]{ A characterization of
varieties whose universal cover is a bounded symmetric domain without ball factors\footnote{AMS Classification: 32Q30, 32N05, 32M15, 32Q20, 32J25, 14C30, 14G35. }}
\author{ Fabrizio  Catanese - Antonio J. Di Scala }

\thanks{The present work took place in the realm of the DFG
Forschergruppe 790 ``Classification of algebraic surfaces and compact
complex
manifolds". The visits of the second author to Bayreuth were supported
by the DFG FOR 790.
The second author was also partially supported by Politecnico di Torino ``Intervento a favore dei Giovani Ricercatori 2009".}

\date{\today}

\begin{abstract}
We give two characterizations of varieties whose universal cover is a bounded symmetric domain without ball factors in terms of the existence of a holomorphic endomorphism $\s$ of the tensor product $T\otimes T^{\vee}$ of the tangent bundle $T$ with the cotangent bundle $T^{\vee}$. To such a curvature type tensor $\s$ one associates the first Mok characteristic cone $\mathcal{CS}$, obtained by projecting on $T$ the intersection of $\mathrm{ker} (\s)$ with the space of rank 1 tensors. The simpler characterization requires that the projective scheme associated to $\mathcal{CS}$ be a finite union of projective varieties of given dimensions and codimensions in their linear spans which must be skew and generate.
\end{abstract}

\maketitle

\section{Introduction}

\noindent
A central problem in the theory of  complex manifolds is  the
one of determining the compact complex manifolds $X$ whose universal
covering
$\tilde{X}$ is biholomorphic to a bounded domain $\Omega\subset \CC^n$.

\noindent
A first important restriction  is given by theorems by   Siegel and
Kodaira (\cite{Kodaira},
\cite{Sie73}) extending to several variables  a  result of Poincar\'e,
and asserting
that necessarily such a manifold $X$ is projective and has ample
canonical divisor $K_X$.

\noindent
A restriction on $\Omega$ is
given by  another  theorem of Siegel (\cite{Siegel}, cf. also
\cite{Kobayashi})  asserting that $\Omega$ must be
holomorphically
convex.

\noindent
The question concerning which domains occur was partly answered by Borel
 ( \cite{Bo63}) who  showed that, given a bounded symmetric domain $\Omega \subset \CC^n$,  there exists a
properly discontinuous group $\Ga \subset \mathrm{Aut}(\Omega)$
   which acts freely on $\Omega$ and is cocompact  (i.e., is such  that
$ X = : \Omega / \Ga$ is a compact complex manifold with
   universal cover $ \cong \Omega$).
\smallskip

\noindent
We consider the following question: given a bounded domain $\Omega \subset \CC^n$, how
can we tell when a projective manifold $X$ with  ample canonical
divisor $K_X$
has $\Omega$ as  universal covering ?

\noindent
The question was solved by Yau (\cite{Yau77})  in the case of a ball, using
the theorem of Aubin and Yau  (see
\cite{Yau78}, \cite{Aubin})  asserting  the existence of K\"ahler Einstein metrics
for varieties with ample canonical bundle. The existence of such  metrics,
joint to some deep knowledge of the differential geometry of bounded symmetric domains,
allows to obtain more general results.

\noindent
Together with Franciosi (\cite{CaFr}) we took up the question for the case of a polydisk,
and a fully satisfactory answer was found in \cite{CaDS} for the special case where
the bounded symmetric domain has all factors of tube type, i.e.,  the domain is
biholomorphic, via the Cayley transform, to some {\bf tube domain} $$ \Omega = V + i  \mathfrak C ,$$
where $V$ is a real vector space and $ \mathfrak C \subset V$ is an open  self dual cone containing no lines.

The main results in the tube case are as follows:

\begin{theo}(\cite{CaDS})
Let $X$ be a compact complex manifold of dimension $n$ with  $K_X$ ample.

Then the
following two conditions (1) and (1'), resp. (2) and (2') are equivalent:

\begin{itemize}

\item[(1)] $X$ admits a slope zero   tensor $0 \neq  \psi   \in
H^0(S^{mn}(\Omega^1_X)(-m K_X) )$, (for some  positive
integer $m$ );

\item[(1')]$X \cong \Omega / \Gamma$ , where
$\Omega$ is a bounded symmetric domain of tube type and $\Gamma$ is a
cocompact
discrete subgroup of
$\mathrm{Aut}(\Omega)$ acting freely.
\item[(2)] $X$ admits a semi special tensor $0 \neq \phi  \in
H^0(S^n(\Omega^1_X)(-K_X) \otimes \eta)$, where $\eta$ is a 2-torsion invertible sheaf,
 such that there is a  point
$p\in X$ for which the
corresponding hypersurface $F_p = :
\{\phi_p = 0 \}\subset \PP (TX_p)$ is  reduced
\item[(2')]
The universal cover of $X$ is a polydisk.
\end{itemize}

Moreover, the degrees and the multiplicities of the irreducible
factors of the polynomial  $\psi_p$ determine uniquely the universal
covering
$\widetilde{X}=\Omega$.

\end{theo}

The main purpose of the present  paper is to extend these results to the more general case
of locally symmetric varieties $X$ whose universal cover is a bounded symmetric domain
without irreducible factors which are isomorphic to a ball of dimension at least two.

\smallskip

In the case where there are no ball factors, we get the following result, using the
concept of an algebraic curvature type
tensor $\s$.

\smallskip
\begin{theo}\label{noball}
Let $X$ be a compact complex manifold of dimension $n$ with $K_X$  ample.

Then the universal covering $\tilde{X}$ is a bounded symmetric domain without factors isomorphic to  higher dimensional balls
if and only if there is a holomorphic tensor $\s \in H^0 ( End (T_X \otimes T_X^{\vee}))$
enjoying the following properties:

1)there is a  point  $p \in X$, and  a splitting of the tangent space $T = T_{X,p}$

$$T = T'_1 \oplus  ... \oplus T'_m$$

such that the first Mok characteristic cone  $\sC\sS$  of $\s$  is $\neq T$ and moreover $\sC\sS$ splits into m
irreducible components $\sC\sS' (j )$  with

2) $ \sC\sS' (j )  =  T'_1 \times  ... \times  \sC\sS'_j \times    ...   \times T'_m$

3)  $\sC\sS'_j  \subset T'_j $  is the cone over a smooth  non-degenerate (that is, the cone $\sC\sS'_j $ spans the vector space $T'_j$) projective variety $\sS'_j$
unless  $ \sC\sS'_j = 0$ and dim $(T'_j) = 1.$

Moreover, we can recover the universal covering  of $\tilde{X}$ from the sequence of pairs
$ ( dim (\sC\sS'_j ) , dim (T'_j))$.

\end{theo}

As we shall recall later,   the first Mok characteristic cone  $\sC\sS \subset  T_{X}$ is defined as the (closure of the) projection on the first factor
of the intersection of $ker (\s )$ with the cone of rank 1 tensors:
$$ ker (\s) \cap \{ t \otimes t^{\vee} \in (T_X \otimes T_X^{\vee}) \}.$$

The above result can be  simplified if we restrict to locally symmetric varieties $X$ whose universal covering
$\tilde{X}$ is a bounded symmetric domain without factors of rank one.

\smallskip
\begin{theo}\label{rank>1}
Let $X$ be a compact complex manifold of dimension $n$ with $K_X$  ample.

Then the universal covering $\tilde{X}$ is a bounded symmetric domain without factors of rank one
if and only if there is   $p \in X$ such that, setting $T = T_{X,p}$,

A) there is a holomorphic tensor $\s \in H^0 ( End (T_X \otimes T_X^{\vee}))$
such that  the first Mok characteristic variety $\sS^1 \subset \PP(T)$ is $ \neq  \PP(T) $ and moreover $\sS^1$  is
the disjoint union of smooth projective varieties
$\sS'_j$ whose projective spans are projectively independent and generate $\PP(T)$.

In other words, iff

A1) there is a  point  $p \in X$ and a splitting of the tangent space $T = T_{X,p}$

$$T = T'_1 \oplus  ... \oplus T'_m$$

such that the first Mok characteristic cone  $\sC\sS$  of $\s$  is $\neq T$ and moreover  $\sC\sS$ splits into m
irreducible components $\sC\sS' (j )$  with

A2) $ \sC\sS' (j )  \subset   T'_j $ and $ \sC\sS' (j )$ generates $  T'_j $

A3)  the projective variety $\sS'_j : = \PP (\sC\sS'_j ) $ is smooth (and non-degenerate, as required in A2)).

Moreover, we can recover the universal covering  of $\tilde{X}$ from the sequence of pairs
$ ( dim (\sC\sS'_j ) , dim (T'_j))$.

\end{theo}

\smallskip
\smallskip

\noindent The above characterizations are important in order to obtain a more precise
formulation of a result of Kazhdan (\cite{Kazh70}).

\noindent
\begin{corollary} Assume that $X$ is a projective manifold with
$K_X$  ample, and that
the universal covering $\tilde{X}$ is a bounded symmetric domain
without irreducible factors which are higher dimensional balls.

   Let $\tau \in
\mathrm{Aut}(\mathbb{C})$ be an automorphism of $\mathbb{C}$.

Then the conjugate variety  $X^{\tau}$ has universal covering
$\tilde{X^{\tau}} \cong \tilde{X}$.
\end{corollary}

It is worthwhile observing that balls of dimension higher than one are taken care of,
once one allows a finite unramified covering,  by  the Yau inequality for summands of the tangent bundle;
hence one can combine the present results with those of
\cite{Yau93} and \cite{V-Z05},  and  obtain  full  results for the general case  where
$\tilde{X}$ is any bounded symmetric domain.
\smallskip

A couple of words about the strategy of the proof:

\begin{enumerate}
\item
knowing that $K_X$ is ample, we have a K\"ahler Einstein metric $h$, and
we consider the Levi-Civita
connection

\item
parallel transport defines then  the restricted
holonomy group $H$ (the connected component of the identity in the holonomy group)
\item
by the theorems of De Rham and Berger (see \cite{Berger} and also
\cite{Olmos}),   the universal covering $\tilde{X}$ of $X$ splits as
 a product
$\tilde{X} = D_1 \times D_2 $ where

\item
$D_1$  is  a bounded symmetric domain without factors of ball type and
$D_2$ is the product of the irreducible factors of dimension $\geq 2$
for which the holonomy group is  the unitary
group;
\item
we observe that by the Bochner principle ( \cite{Koba80}) the tensor $\s$ is parallel,
hence if we restrict the tensor at any point $p\in X$ we observe that $\s_p$
is $H$-invariant.
\item
We decompose the holonomy group as a product and accordingly the vector space $T$, the tangent space
to $X$ at $p$.
\item
We use elementary representation theory to derive some restrictions which the tensor
$\s$ must satisfy (this is done in section 3)
\item
we associate to any such tensor $\s$ its first Mok characteristic variety (see section 2), having in mind two standard
possibilities constructed using the algebraic curvature tensors of irreducible
bounded domains  defined by Kobayashi and Ochiai
in \cite{KobOchi}.
\item
the rest of the proof (sections 4 and 5) is projective geometry, using Mok's description
(\cite{Mok1}, \cite{MokLibro}) of the orbits on $T_i \otimes T_i^{\vee}$
of the complexified holonomy group of an irreducible bounded domain.

\end{enumerate}

\section{Algebraic curvature-type tensors and their First Mok characteristic varieties}

In this section we consider the following situation.
We are given a direct sum

$$T = T_1 \oplus  ... \oplus T_k$$
of irreducible representations $T_i$ of a group $H_i$
(the unusual notation is due to the fact that $T$ in the application shall be the tangent space to a projective manifold
at one point, and $H = H_1 \times \dots \times H_k$ shall be the restricted holonomy group).

\begin{defin}
1)
An algebraic curvature-type tensor is a nonzero  element
$$ \s   \in End (T \otimes T^{\vee}).$$

2) Its first Mok characteristic cone  $\sC\sS \subset  T$ is defined as  the projection on the first factor
of the intersection of $ker (\s)$ with the set of rank 1 tensors, plus the origin:

$$ \sC\sS : = \{ t \in T | \exists t^{\vee} \in T^{\vee} \setminus \{0\}, (t \otimes  t^{\vee} ) \in ker (\s)  \}.$$

3) Its {\bf first Mok characteristic variety} is the subset $\sS : = \PP ( \sC\sS) \subset \PP (T)$.

4) More generally, for each integer $h$, consider
$$ \{  A \in T \otimes T^{\vee} | A  \in ker (\s), {\rm Rank} (A) \leq h  \},$$
and consider the algebraic cone which is its  projection on the first factor
$$ \sC\sS^h : = \{ t \in T | \exists A  \in ker (\s), {\rm Rank} (A) \leq h  , \exists t'  \in T: t = A   t'  \},$$
and define then $\sS^h : = \PP ( \sC\sS^h) \subset \PP (T)$ to be the {\bf h-th Mok characteristic variety.}

5) We define then the {\bf full characteristic sequence} as the sequence
$$ \sS = \sS^1   \subset \sS^2 \subset \dots \subset \sS^{k-1} \subset \sS^k =  \PP (T).$$
\end{defin}

\begin{remarkk}\label{Mok}
\begin{enumerate}
\item
In the case where $\s$ is the curvature tensor of an irreducible symmetric bounded domain $\sD$, Mok (\cite{Mok1})
proved that the difference sets  $ \sS^h \setminus    \sS^{h-1}$ are exactly all the  orbits of the  parabolic subgroup $P$
 associated to the compact dual
$\sD^{\vee} = G/P$.  In particular, the algebraic cone  $ \sC \sS^h$ is irreducible and $H_i$-invariant.
\item
More generally, if $ \sS$ is an $H_i$-invariant algebraic cone, and  $H_i$ is the holonomy group of an  irreducible symmetric bounded domain,
then necessarily $ \sS$ is irreducible and indeed equal to one of the $ \sS^h$.
\item
In the case instead where $H_i$ acts as the full unitary group on $T_i$, then any $H_i$-invariant algebraic cone in $T_i$ is trivial, that is, either equal to
$T_i$ or just equal to  $\{ 0 \}$, where $0 \in T_i$ is the origin.

\end{enumerate}

\end{remarkk}

\begin{lemma}\label{MokSmooth}
If  $\s$ is the curvature tensor of an irreducible symmetric bounded domain $\sD$,   $ \sS^h$ is smooth if and only if $h=1$.

\end{lemma}

\Proof
That $ \sS^1$ is smooth follows form the above remark, since $ \sS^1$ is a single orbit.

Conversely, observe that we have a sequence of inclusions for the Mok characteristic varieties:
$$ \sS = \sS^1   \subset \sS^2 \subset \dots \subset \sS^{k-1} \subset \sS^k =  \PP (T).$$

Let then $P \in  \sS^1$, and let $G$ be the stabilizer of $P$. The tangent space $V$ of $ \PP (T)$ at $P$
is a $G$-representation, and the Zariski tangent spaces to $ \sS^h$ yield a flag of $G$-invariant subspaces of $V$,
$$0 \subset V_1 \subset V_2 \subset \dots \subset V.$$

Our assertion follows then from the claim that
$V_1$ and $ V$ are the unique invariant subspaces;
 since then, for   $k >  h > 1$, we have
$ V_h \neq V_1 \Rightarrow V_j = V$ therefore $ \sS^h$ is singular at every point of $ \sS^1$.

Let us prove the claim.

The first characteristic variety $ \sS^1$ is homogeneous for the compact subgroup $K$ which
stabilizes the origin of the bounded domain $\sD$.
The stabilizer $ K_P$, as proven   in Theorem 2.3 and Proposition 2.4 at page 5 of \cite{CoDS}, acts irreducibly on the normal space
to $ \sS^1$ at $P$.  Therefore $V$ splits as $V_1 \oplus N_P$, where $N_P$ is
$ K_P$- irreducibile.  Now, $G$ contains $ K_P$, hence  if  $W$  is $G$ -invariant and strictly contains  $ V_1$,
then   $W$ is also $K_P$- invariant and   $W = V$.

\qed

One geometric situation we have particularly in mind is the one where
$$ (**) \  \s  = \oplus_{i=1}^k \s_i  \in  \oplus_{i=1}^k (End (T_i \otimes T_i^{\vee})) \subset End (T \otimes T^{\vee}).$$

In this case it is clear that

$$  ker (\s) =   \oplus_{i=1}^k  ker (\s_i) \bigoplus ( \oplus_{i \neq j} (T_i \otimes  T_j^{\vee})),$$

and if we intersect the Kernel of $\s$ with the set of rank 1 tensors, we obtain

$$ ker (\s) \cap \{ t \otimes  t^{\vee} \} = \{  t \otimes  t^{\vee}| t = \Si_{i=1}^k t_i ,  t ^{\vee}= \Si_{j=1}^k t^{\vee}_j ,  \forall i \
( t_i \otimes  t_i^{\vee}) \in ker (\s_i)\}.$$

Defining now

$$ \sC\sS_i : = \{ t_i \in T_i | \exists t_i^{\vee} \in T_i^{\vee} \setminus \{0\}, (t_i \otimes  t_i^{\vee} ) \in ker (\s_i)  \},$$
and similarly $\sS_i : = \PP ( \sC\sS_i) \subset \PP (T_i)$
we see therefore that under hypothesis (**) we have an inclusion
$   \oplus_{i=1}^k  \sC\sS_i  \subset  \sC\sS  $.

But  indeed, since  $$  \sC\sS = \{ t = \Si_j t_j | \exists t^{\vee} \neq 0 , t_i \otimes t_i^{\vee} \in ker (\s_i) \forall i \}=
 \{ t = \Si_j t_j | \exists t_i^{\vee} \neq 0 , t_i \otimes t_i^{\vee} \in ker (\s_i) \},$$
 we have that
$$  \sC\sS  =  \cup _{i=1}^k  T_1 \oplus T_2 \dots \oplus T_{i-1} \oplus \sC\sS_i  \oplus T_{i+1} \oplus \dots
\oplus T_k .$$

The above formula yields a decomposition of the Zariski closed projective set $  \sS $ as the union of the
Zariski closed projective sets

$$  \sS(i) :   =  \PP (T_1 \oplus T_2 \dots \oplus T_{i-1} \oplus \sC\sS_i  \oplus T_{i+1} \oplus \dots
\oplus T_k) .$$

The latter sets are the join of the linear subspace $$\PP(T_{\hat{i}}) : =  \PP (T_1 \oplus T_2 \dots \oplus T_{i-1}  \oplus T_{i+1} \oplus \dots
\oplus T_k) $$ with $\sS_i$.

\begin{rem}
The next question is: when is the above an irredundant decomposition?

It is a necessary condition that each $\sC\sS_i  \neq T_i$ (i.e., $\sS_i \neq \PP (T_i)$),
otherwise $\sC\sS (i) = T$.

This condition is also sufficient. In  fact,   the irreducible components of $\sS (i) $ are the joins of $\PP(T_{\hat{i}}) $
with the irreducible components of $\sS_i$, hence for each  component  of $\sS (i) $ the projection onto $\PP(T_j)$  is surjective
whenever $j \neq i$: therefore this component cannot be contained in any $\sS (j) $ when $ j \neq i$.

\end{rem}

\begin{rem}
At the other extreme, if $\s_i = Id_{T_i}$, or $\s_i$ is invertible, then $ ker (\s_i )= 0$, hence in this case
$\sS (i) = \PP(T_{\hat{i}}) $.

More generally, it can happen that $\sS (i)$ is a linear subspace, iff  $\sC \sS_i$ is a linear subspace.
In the sequel, we shall assume that each $\s_i$ is $H_i$-invariant: hence $\sC \sS_i$ shall be an invariant subspace
of $T_i$: by the irreducibility of $T_i$, the only possibility is either that $\sC \sS_i = 0$, or that  $\sC \sS_i = T_i$.

We can avoid both possibilities by requiring (the second case should only occur for factors $T_i$ of dimension = 1)

\begin{enumerate}
\item
$\sS \neq \PP(T)$
\item
if $\sS$ has a component $\sS^0$ which is a linear subspace, then $\sS^0$ must be a hyperplane.

\end{enumerate}

\end{rem}

\section{Holonomy invariant curvature-type tensors}

In this section we continue our consideration of a curvature-type tensor $\s $ assuming that it is invariant by the natural action
of the group $H = H_1 \times \dots \times H_k$.

We can naturally write $\s$ as a direct sum
$ \s = \oplus_{(i,j),(h,k)} \s_{(i,j),(h,k)}$,
$$ \s_{(i,j),(h,k)} \colon T_i \otimes  T_j^{\vee} \ra T_h \otimes  T_k^{\vee}.$$

\begin{lemma}
$ \s_{(i,j),(h,k)} = 0 $ if $i \neq h$ or $j\neq k$, while for $ i\neq j$
$$ \s_{(i,j)}: =  \s_{(i,j),(i,j)} \colon T_i \otimes  T_j^{\vee} \ra T_i \otimes  T_j^{\vee}$$
is a multiple of the identity.
\end{lemma}

\Proof
The second assertion is a consequence of Schur's lemma once we show that $T_i \otimes  T_j^{\vee} $
is, for $ i \neq j$, an irreducible representation of the compact group $H_i \times H_j$.

We use moreover that $T_i$ is an irreducible representation of $H_i$.
$H_i$ being compact, if $\chi_i$ is the character of the representation $T_i$, and $d\mu_i$ is the Haar measure
of $H_i$, we have that irreducibility
is equivalent to
$$ \int _{H_i} \chi_i \overline{ \chi_i} \  d\mu_i = 1
$$
Since the character $\chi_{i,j}$ of $T_i \otimes  T_j^{\vee} $ on $H_i \times H_j$ is
$$  \chi_{i,j} (x,y) := \chi_{i}(x)\overline{ \chi_{j} (y)}$$
and
$$ \int _{H_i \times H_j }  |\chi_{i,j} (x,y)|^2  d\mu_{i,j} = {\rm [by \ Fubini]} \ =  \int _{H_i} \chi_i \overline{ \chi_i} \  d\mu_i
\cdot  \int _{H_j} \overline{ \chi_j} \chi_j \  d\mu_j  = 1$$
we conclude that $T_i \otimes  T_j^{\vee} $
is, for $ i \neq j$, an irreducible representation of  $H_i \times H_j$.

For the first assertion, assume now that $i \neq h$ and let
$$\s' : =   \s_{(i,j),(h,k)} \colon T_i \otimes  T_j^{\vee} \ra T_h \otimes  T_k^{\vee}.$$

Since $\s$ is $H$-invariant, $\s'$ is $H_i \times H_h$-invariant. By what we have seen
$Hom (T_i, T_h) \cong T_h \otimes  T_i^{\vee} $ is an irreducible nontrivial representation
of $H_i \times H_h$, hence there are no $H_i \times H_h$ invariant homomorphisms in $Hom (T_i, T_h) $.
A fortiori $\s' = 0$.

A completely analogous argument yields  $\s' = 0$ if $j \neq k$.

\qed

Using the previous lemma we consider the first characteristic variety in the case where
$\s$ is $H$-invariant.

In this case it is clear that

$$  ker (\s) =   \oplus_{i=1}^k  ker (\s_i) \bigoplus ( \oplus_{i \neq j, \s_{i,j} = 0} (T_i \otimes  T_j^{\vee})),$$

and if we intersect the Kernel of $\s$ with the set of rank 1 tensors, we obtain

$$ ker (\s) \cap \{ t \otimes  t^{\vee} \} =$$
$$= \{  t \otimes  t^{\vee}| t = \Si_{i=1}^k t_i ,  t ^{\vee}= \Si_{j=1}^k t^{\vee}_j ,  \forall i \
( t_i \otimes  t_i^{\vee}) \in ker (\s_i), t_i \otimes t^{\vee}_j=0 \  \forall \ i \neq j \  s.t. \  \s_{i,j} \neq 0 \}.$$

Its projection on the first factor is the set
{\small
$$\sC\sS = \cup_j  \sC\sS (j)  : = \cup_j  \{  t = \Si_{i=1}^k t_i   \in T | \exists t_j^{\vee} \in T_j^{\vee} \setminus \{0\}, (t_j \otimes  t_j^{\vee} ) \in ker (\s_j),
t_i = 0   \forall \ i \neq j \  s.t. \  \s_{i,j} \neq 0  \} ,$$}
and we have a corresponding set  $\sS (j)  : = \PP ( \sC\sS (j) ) \subset \PP (T)$.

We can also write

$$\sC\sS (j)  =  \hat{T}_{1,j} \oplus \hat{T}_{2,j} \dots \oplus \hat{T}_{j-1,j} \oplus \sC\sS_j  \oplus \hat{T}_{j+1,j} \oplus \dots
\oplus \hat{T}_{k,j}  ,$$
where $\hat{T}_{i,j} = T_i $ if $\s_{i,j} =  0$, $\hat{T}_{i,j} = 0$ if $\s_{i,j} \neq 0$.
In other words, if one forgets about the ordering,
$$\sC\sS (j)  =  \sC\sS_j  \bigoplus (\oplus_{\s_{i,j} = 0} T_i  ).$$

Our first remark is that, in case where there exists a $ \s_{i,j} \neq 0$ with $i\neq j$, then necessarily $\sS (j)$
is contained in a hyperplane.

Our second observation is that however in this case the decomposition $\sC\sS = \cup_j  \sC\sS (j) $
does not need to be irredundant, as one sees already in the case $k=2$, $ \s_{1,2} \neq 0$, $ \s_{2,1} = 0$.

We end this section with an important, even if trivial, remark.

\begin{rem}\label{invsubspace}
If $V$ is an $H$-invariant linear subspace of $T$, then there is a subset $ \sI \subset \{ 1, \dots , k\}$ such that
$ V = \oplus_{ i \in \sI} T_i$.
\end{rem}

\section{Proof of Theorem \ref{noball}}

One implication follows right away from section 3, since we may take $\s : = \oplus_{i=1}^k \s_i$, letting
$\s_i$  be  the algebraic curvature tensor of a bounded symmetric domain of rank  greater than one (cf. \cite{KobOchi},
lemma 2.9), and the identity of $T_i \otimes T_i^{\vee}$ for the factors of dimenson equal to one.
Then $$  \sC\sS  =  \cup _{i=1}^k  ( T_1 \oplus T_2 \dots \oplus T_{i-1} \oplus \sC\sS_i  \oplus T_{i+1} \oplus \dots
\oplus T_k ).$$

The converse implication follows since $K_X$ is ample, hence we may consider the K\"ahler-Einstein metric of $X$,
for which the tensor $\s$ is parallel, as proven by Kobayashi in (\cite{Koba80}) (since it is a tensor of covariant type
two and contravariant type also two).

Hence the (irredundant)  irreducible decomposition

$$\sC\sS = \cup_{j=1}^m \sC\sS' (j )  = \cup_{j=1}^m ( T'_1 \times \dots \times   \sC\sS'_j \times    ...   \times T'_m)$$
is invariant under the holonomy group $H$.

 Observe that, since $K_X$ is ample, all the irreducible holonomy factors $H_i$ are either equal to $ U (T_i)$,
or $H_i$ acts on $T_i$ as the holonomy of a bounded symmetric domain. This implies that
if  $\sC\sS_j \subset T_j $ is a proper $H_j$-invariant subset, then $\sC\sS_j $ is (see, remark  \ref{Mok}, and also
 \cite{Mok1}, and also \cite{MokLibro}, page 252) a Mok characteristic variety and $T_j$ is a bounded symmetric domain factor.

The holonomy invariance implies, as shown in section 4, that there is another (possibly redundant) irreducible
(by what we have just observed) decomposition
$$\sC\sS = \cup_{i=1}^k \sC\sS (i)    = \cup_{i=1}^k ( \sC\sS_i  \bigoplus (\oplus_{\s_ {j,i} = 0} T_j  )).$$

It follows immediately that $ k \geq m$.

On the other hand, by our assumption and by Lemma \ref{smooth}, the linear subspace
$$ \tilde{T'}_j : = ( T'_1 \oplus \dots \oplus T '_{j-1} \oplus  T '_{j+1}  \oplus  \dots   \oplus T'_m)$$
is the maximal vector subspace $V$ such that  $V + \sC\sS' (j ) \subset \sC\sS' (j )$, hence these linear subspaces
are holonomy invariant, in particular their mutual intersections are holonomy invariant.

We conclude that each subspace $T'_j$ is holonomy invariant.
By Remark \ref {invsubspace} each $T'_j$  is a sum of a certain number of $T_i $'s.

Comparing the two decompositions, it follows that each $\sC\sS' (j ) $ equals some $ \sC\sS ( i  ) $, and the hypothesis that the
linear span of $\sC\sS' (j ) $ equals $T$ implies that
$$\sC\sS' (j )  =  \sC\sS ( i  ) =  \sC\sS_i  \bigoplus (\oplus_{j \neq i} T_j  ) = : \sC\sS_i  \oplus \tilde{T}_i.$$

Once more, by Lemma \ref{smooth} $ \tilde{T'}_j $ is the maximal linear subspace $V$ such that $ V + \sC\sS' (j )  \subset \sC\sS' (j ) $,
and the above equality show that this subspace contains  $ \tilde{T}_i.$
Since all the subspaces $T'_j$ yield a direct sum,and are holonomy invariant,  it follows that  $T'_j =T_i$, and $\sC \sS'_j = \sC \sS_i $.

Therefore we finally obtain that $m =k$ and that,  when $\sC \sS'_j = \sC \sS_i  \neq 0$,
then $ \sS'_j =  \sS_i $
 is a smooth projective variety.

Since  the only smooth characteristic variety is the first Mok characteristic variety (as shown in Lemma \ref{MokSmooth}), it follows that the cones
$ \sC \sS_i $ are just the origin when $ \dim (T_i) = 1$, or they are the cones over the first Mok characteristic variety.

To finish the proof,  we must only show the following claim

{\bf Claim}\label{Mok}  The dimension and codimension of the  first Mok characteristic variety determines the
irreducible bounded symmetric domain $\sD$  of rank $\geq 2$.

\Proof This follows from the following table.

Let $\sD$ be an irreducible Hermitian symmetric space of rank $> 1$.

The following table follows from Mok's enumeration of the characteristic variety $\sS^1(\sD)$, see Mok's Book \cite{MokLibro}, page 250.

\vspace{1cm}

 \begin{tabular}{|c|c|c|}
  \hline
  $\sD$ & $\mathrm{dim}(\sD)$ & $\mathrm{dim}(\sS^1(\sD))$ \\
  \hline
  $I_{p,q}$ & $pq$ &  $p + q - 2$ \\
  \hline
  $II_n$ & $\frac{n(n-1)}{2}$ & $2(n-2)$\\
  \hline
  $III_n$ & $\frac{n(n+1)}{2}$ & $n-1$\\
  \hline
  $IV_n$ & $n$ & $n-2$\\
  \hline
  $V$ & $16$ & $10$\\
  \hline
  $VI$ & $27$ & $16$\\
  \hline
\end{tabular}

\vspace{1cm}

Let $\eta : \mathcal{IHSS} \to \mathbb{N} \times \mathbb{N}$ be defined as \[ \eta(\sD) := (\mathrm{dim}(\sD) , \mathrm{dim}(\sS^1(\sD))) \]

{\bf Fact:} $\eta$ is injective.\\

\Proof The proof is obtained by direct inspection of the above table. Indeed, it is not difficult to check that the pairs $(27,16)$ and $(16,10)$ comes just from the domains of $VI$ and $V$ respectively.
To show that the pair comings from domains $IV_n$ come just from the domains of type $IV_n$ it is necessary to recall the following isomorphisms:
\[ IV_3 \cong III_2 \, , \,  IV_6 \cong II_4 \, , \,  IV_4 \cong I_{2,2} \]

\qed

\section{Proof of Theorem \ref{rank>1}}

In one direction, if $X$ is locally symmetric without factors of rank 1,
consider the tensor $\s$ such that $\s_i$ is the curvature tensor for all $i$,
and $\s_{i,j} $ is the identity on $T_i \otimes T_j^{\vee}$  $\forall \ i \neq j$.

We saw then in section 4 that
$$\sC\sS (j)  =  \sC\sS_j  \bigoplus (0)$$

and then $\sS (j) \subset \PP(T_j)$ is the first Mok characteristic variety, which is smooth, hence A1), A2) and A3) hold.

Conversely, all the components $\sC\sS (j)$ of the cone $\sC\sS$ contain no nontrivial vector subspace, since the cone
over a projective variety is singular unless the variety is a linear subspace.
Hence, by the observations made in section 4 it follows that the holonomy invariant tensor $\s$
is such that all the components  $\s_{i,j} $ are a nonzero multiple of  the identity on $T_i \otimes T_j^{\vee}$  $\forall \ i \neq j$.

Then $\sC\sS (j)  =  \sC\sS_j  \bigoplus (0)$ and the projective variety $\sS_j $ is smooth and holonomy invariant, therefore
we conclude as for theorem \ref{noball} that $\sS_j $ is the first Mok characteristic variety, and that we recover
the universal cover from the variety $\sS$.

\section{Proof of the Kazhdan's type corollary}

Consider the conjugate variety $X^{\tau}$: since $K_X$ is ample we
may assume that $X$ is projectively embedded by $H^0 ( X, \hol_X (m
K_X)$.

$\tau$ carries $X$ to
$X^{\tau}$ and $K_X$ to $K_{ X^{\tau}}$, hence also $X^{\tau}$
has ample canonical divisor.

Moreover, $\tau$ carries the algebraic curvature type tensor $\s$ to a similar tensor
$\s^{\tau}$. The equations of the Mok characteristic varieties are defined over $\ZZ$,
hence we obtain that $\tau$ transforms each variety $\sS^i(X, \s)$ into $\sS^i(X^{\tau}, \s^{\tau})$,
in particular respecting their dimension and codimension.

We conclude then immediately by the last assertion of our main
theorems that  the
universal covering of
$X^{\tau}$ is $\tilde{X}$.

\qed

\section{Elementary lemmas}

We collect here, for the readers' benefit, some trivial but  important observations.

\begin{lemma}\label{smooth}
Let $\sS \subset \PP (V) = \PP^n$ be a non-degenerate projective variety, $ \sS \neq \PP(V)$,
and consider the join $Z : =  \sS * \PP (W) \subset \PP (V \oplus W) = \PP^{n+m}$.
Then   $Z$ is smooth $\Leftrightarrow$ $W = 0$ and $\sS$ is smooth.
\end{lemma}

\Proof
Let $I$ be the homogeneous ideal of $\sS$. Since $\sS$ is non-degenerate, each $ 0 \neq f \in I$
has degree $\geq 2$, and moreover $I $ contains some non zero polynomial (since $ \sS \neq \PP(V)$).

We shall show that $\PP(W) \subset Sing (Z)$, observing that $\PP(W) \neq \emptyset $ unless $ W=0$, which is exactly which we have to proof.

Observe that

$$ Z = \{ (v,w) | f (v) = 0, \forall f \in I \}. $$
Hence

$$ Sing (Z) = \{ (v,w) | \frac{\partial f }{\partial v_j} (v) = 0, \forall f \in I \}. $$
Since however $ deg (f) \geq 2 $ , $\forall f \in I, f \neq 0$, $ \frac{\partial f }{\partial v_j} (v)$ vanishes for $ v=0$,
hence  $\PP(W)  = \{ (0,w) | w \in W \} \subset Sing (Z)$.

\qed

In the next lemma we use our standard notation, introduced in section 2.

\begin{lemma}\label{subspace}
Let $\sC \sS \subset  T_i $ be an $H_i$ -invariant algebraic cone.

Then there is no nontrivial linear subspace $V_i$ such that  $ V_i + \sC \sS \subset \sC \sS $,
unless $\sC \sS  =   T_i $.

\end{lemma}

\Proof
If $V_i$ is nontrivial, then $W_i : = \{ v | v  + \sC \sS \subset \sC \sS \}$ is a non trivial linear subspace, which is
$H_i$ invariant. But $T_i$ is an irreducible representation, hence if $W_i \neq \{ 0\}$, then  $\sC \sS  =   T_i $.

\qed

{\bf Acknowledgement:}

We would like to heartily thank Ngaiming Mok  for drawing our attention to  the reference \cite{KobOchi}
and putting us on the right track for the generalization of the results of \cite{CaDS}.

 Thanks also to Daniel Mckenzie for pointing out that our first proof of  Theorem  \ref{noball}  was erroneously written.

\medskip
\noindent {\bf Authors '  Addresses:}\\
\noindent Fabrizio Catanese, \\ Lehrstuhl Mathematik VIII, Mathematisches
Institut der \\Universit\"at Bayreuth\\ NW II,  Universit\"atsstr. 30,
95447
Bayreuth, Germany.\\
            email:          fabrizio.catanese@uni-bayreuth.de\\

\noindent Antonio J. Di Scala, \\ Dipartimento di Scienze Matematiche ``G.L. Lagrange",\\
Politecnico di Torino, \\ Corso Duca degli Abruzzi 24, 10129 Torino,
Italy. \\
    email:     antonio.discala@polito.it \\

\end{document}